\documentclass[12pt]{amsart}
\usepackage{amssymb}
\usepackage[all]{xy}

\newcommand{\issuenumber}{13}
\newcommand{\issuemonth}{June}
\newcommand{\issueyear}{2005}

\newtheorem{issue}{Issue}

\theoremstyle{definition}

\theoremstyle{remark}

\newcommand{\ed}{\end{thebibliography}\general\end{document}}

\renewcommand{\b}{\mathfrak{b}}

\newcommand{\p}{\mathfrak{p}}

\newcommand{\NON}{{\mathsf   {NON}}}
\newcommand{\COF}{{\mathsf   {COF}}}

\newcommand{\M}{\mathcal{M}}

\newcommand{\cov}{\mathsf{cov}}

\newcommand{\fo}{\mathfrak{od}}

\renewcommand{\b}{\mathfrak{b}}

\renewcommand{\split}{\mathsf{Split}}
\newcommand{\bq}{\begin{quote}}
\newcommand{\eq}{\end{quote}}

\newcommand{\B}{\mathcal{B}}
\newcommand{\BG}{\B_\Gamma}

\newcommand{\sone}{\mathsf{S}_1}    \newcommand{\sfin}{\mathsf{S}_{fin}}
\newcommand{\ufin}{\mathsf{U}_{fin}}

\newcommand{\nin}{\not\in}

\newcommand{\naturals}{{\mathbb N}}
\newcommand{\N}{\naturals}
\newcommand{\sm}{\setminus}

\newcommand{\by}[2]{\par\hfill\emph{#1}, #2}
\newcommand{\nby}[1]{\par\hfill\emph{#1}}
\newcommand{\Tau}{\mathrm{T}}
\newcommand{\CE}{\textsc{CE}}

\newcommand{\be}{\begin{enumerate}}
\newcommand{\ee}{\end{enumerate}}
\newcommand{\bi}{\begin{itemize}}
\newcommand{\ei}{\end{itemize}}
\renewcommand{\i}{\item}

\newcommand{\general}{\small\vfill\par\noindent\hrulefill\par
\noindent\textbf{Previous issues.} The first issues of this
bulletin, which contain general information (first issue), basic
definitions, research announcements, and open problems (all
issues) are available online, on \arx{math.GN/$x$}, where $x$ is
\texttt{0301011}, \texttt{0302062}, \texttt{0303057},
\texttt{0304087}, \texttt{0305367}, \texttt{0312140},
\texttt{0401155}, \texttt{0403369}, \texttt{0406411},
\texttt{0409072}, \texttt{0412305}, and \texttt{0503631},
respectively, for issues number $1$ to $12$.\\[0.1cm]
\textbf{Contributions.}
Please submit your contributions (announcements, discussions, and open problems)
by e-mailing us. It is preferred to write them
in \LaTeX{}.
The authors are urged to use as standard notation as possible, or otherwise give
the definitions or a reference to where the notation is explained.
Contributions to this bulletin would not require any transfer of copyright,
and material presented here can be published elsewhere.\\[0.1cm]
\textbf{Subscription.}
To receive this bulletin (free) to your
e-mailbox, e-mail us:\\
{boaz.tsaban@weizmann.ac.il}
}

\newcommand{\nArxPaper}[5]{\subsection{#2}{#4}\par\hfill{\arx{#1}}\par\hfill\emph{#3}}

\newcommand{\nAMSPaper}[5]{\subsection{#2}{#4}\par\hfill{\texttt{#1}}\par\hfill\emph{#3}}
\newcommand{\SPMBul}{\textbf{$\mathcal{SPM}$ Bulletin}}

\newcommand{\arx}[1]{\texttt{http://arxiv.org/abs/#1}}
\newcommand{\url}[1]{\bq\texttt{#1}\eq}
\newcommand{\online}[1]{The paper is available online at \url{#1}}

\title[$\mathcal{SPM}$ Bulletin \textbf{\issuenumber} (\issuemonth{} \issueyear)]{%
$\mathcal{SPM}$ Bulletin\\[0.5cm]
Issue number \issuenumber: \issuemonth{} \issueyear{} \CE{}}

\begin{document}
\maketitle

\tableofcontents

\section{Editor's note}

Please mark your calendars: 19--22 December 2005 are the days of the coming
large workshop on Selection Principles in Mathematics.
The current list of participants (to be available in the Workshop's homepage soon)
promises that this will be an exceptionally interesting workshop.
We hope to see many of the readers of this bulletin there.
See Section \ref{spmc}.

In Section \ref{cichon}
we announce another, very close (by time and theme), fascinating workshop,
to be held in 24--31 July 2005.

This issue also contains, as usual, research announcements.

\medskip

Contributions to the next issue are, as always, welcome.

\medskip

\by{Boaz Tsaban}{boaz.tsaban@weizmann.ac.il}

\hfill \texttt{http://www.cs.biu.ac.il/\~{}tsaban}

\section{Workshops on SPM themes}

\subsection{Second workshop on Coverings, Selections and Games in Topology (SPM05)}\label{spmc}
Lecce, Italy, 19--22 December 2005.

The study of Selection Principles in Mathematics has experienced
rapid expansion during the past few years with a large number of
mathematicians contributing to the area, and entering the area.
The combination of classical and modern methods has lead to
fascinating breakthroughs and to complete solutions of some of the
oldest open problems (1920's and 1930's) in the field, and a large
number of new problems covering a variety of topics in mathematics
have been identified. Since 2001 several strong Ph.D.\ theses were
devoted exclusively to topics in Selection Principles in
Mathematics.

Though Selection Principles in Mathematics had its original
beginnings mostly in the study of covering properties of
topological spaces that were introduced by Menger (1924), Hurewicz
(1925), Rothberger (1937) and Sierpinski (1937), the field has
become vastly wider. There are currently several well-defined
focus areas in Selection Principles in Mathematics, including:
\be
\i Distributivity properties in Boolean algebras
\i Combinatorial properties of filters on the natural numbers
\i Boundedness properties in topological groups
\i Closure- and convergence- properties in function spaces
\i Combinatorial cardinal characteristics of the continuum
\i Selective screenability and covering dimension
\i Covering properties of topological spaces.
\ee
The aim of the workshop is to survey current directions in the
field through a number of plenary talks and to learn about current
results and open problems in this area through a number of shorter
contributed talks.

A tentative list of plenary speakers includes:
\begin{itemize}
\i Liljana Babinkostova (Boise State University)
\i Taras Banakh (Lviv University)
\i Lev Bukovsky (P.\ J.\ Safarik University)
\i Filippo Camorroto (University of Messina)
\i Ljubisa Kocinac (University of Nis)
\i Giuseppe Di Maio (Second University of Napoli)
\i Heike Mildenberger (Kurt G\"odel Research Center for Mathematical Logic)
\i Arnold Miller (University of Wisconsin)
\i Masami Sakai (Kanagawa University)
\i Marion Scheepers (Boise State University)
\i Boaz Tsaban (Weizmann Institute of Science)
\i Lubomyr Zdomsky (Lviv University)
\end{itemize}

The registration fee for the conference is $50$ Euros. There is
limited financial support available for participation. It could
cover the registration fee, accommodation or meals. People who
have no other financial support for participating in the
conference should apply \textbf{before the end of June, 2005}, to the
Organizing Committee for financial support. Direct all
applications to Professor Cosimo Guido (\texttt{cosimo.guido@unile.it}).

Information on travelling to Lecce and on accommodations will
become available in the near future at the workshop web sites
currently under construction. The url's are:
\url{http://www.matematica.unile.it/mostra\_avviso.asp?n=136}
\url{http://diamond.boisestate.edu/\~{}spm/Lecce2/index.htm}
Topology Atlas is generously providing abstract services for this
conference. Please submit abstracts at the following web-site:
\url{http://atlas-conferences.com/cgi-bin/abstract/submit/caqh-01}
Please note that the deadline for submitting abstracts is
\textbf{October 31, 2005}.

Submitted abstracts can be viewed at
\url{http://atlas-conferences.com/cgi-bin/abstract/caqh-01}
The current list of sponsors for the meeting includes:
The University of Lecce,
Department of Mathematics ``E. De Giorgi'' -- University of Lecce, and
Topology Atlas.

\nby{L.\ Babinkostova,
C.\ Guido,
L.\ Kocinac,
M.\ Scheepers, and
B.\ Tsaban}

\subsection{Analysis and Descriptive Set Theory Workshop}\label{cichon}
24--31 July 2005, Banach Center, 
Bedlewo,
 Poland.

\subsubsection*{Descriptive set theory: Effective methods, equivalence relations}
Speakers: J.\ Cicho\'n (Wroc{\l}aw, Poland) and S.\ Solecki (Urbana-Champaign,USA).

Recent applications of descriptive set theory emphasize
investigation of definable equivalence relations in various
mathematical contexts (as opposed to the study of definable sets)
and usage of effective (recursive theoretic) methods. At the root
of this development lie, on the one hand, the theorem of Silver
that each co-analytic equivalence relation on a Polish space
either has countably many classes or there exists a perfect set of
inequivalent elements and, on the other hand, the effective
descriptive set theoretic proof of this theorem due to Harrington.
Extensions of Silver's theorem and of its method of proof have
found numerous applications in Polish group actions, in certain
classification problems in algebra and topology, in continuum
theory, and in real function theory, to name only a few.

In the workshop, we propose to present Harrington's proof of
Silver's theorem along with a survey of a broader mathematical
background of this result. The talks will be aimed at graduate
students interested in applications of descriptive set theory.
Prior exposure to basic facts of classical descriptive set theory
will make motivation of some of the material clearer.

Here is the plan of the first half of the workshop.

\begin{enumerate}
\item Rudiments of recursion theory through Kleene's recursion
theorem.

\item Survey of basic theorems of classical descriptive set theory
(Borel, analytic, co-analytic sets, operation $\mathcal A$,
Nikodym's theorem on preservation of Baire property by operation
$\mathcal A$, Mycielski's theorem on perfect independent sets for
meager relations).

\item Kleene's classes $\Sigma^0_1$, $\Pi^0_1$, $\Sigma^1_1$,
$\Pi^1_1$, good universal sets.

\item $\Pi^1_1$ sets and their representation via well-orderings,
boundedness, reflection theorem, coding $\Delta^1_1$ sets, Gandy's
basis theorem for $\Sigma^1_1$ sets.

\item The Gandy-Harrington topology and the proof of Silver's
theorem.
\end{enumerate}

\subsubsection*{Analysis: typical functions, level sets structure}

Speakers: U.\ B.\ Darji (Louisville, USA) and M.\ Morayne (Wroc{\l}aw, Poland)

The purpose of this course is to study properties of smooth
functions in terms of their level sets. We will first recall the
Baire Category Theorem and the celebrated Banach Theorem that a
typical continuous function on $[0,1]$ is nowhere differentiable.
Using the Baire Category theorem, we will also prove the existence
and the abundance of otherwise complicated objects in various
branches of analysis and topology. However, the main thrust is to
study the behavior of a typical smooth function in terms of its
level sets (inverse images of points) structure. This topic was
opened by a well known theorem of Bruckner and Garg describing the
level sets of typical continuous functions defined on $I$. The
proof of this theorem will be presented as well as the proofs of
recent results concerning level sets of typical $C^n$ functions,
$1 \le n \le \infty$. We will also discuss level sets structure of
a typical continuous function from $S^2$, the 2-sphere, into $I$.
This is more of a topological result which shows that objects such
as pseudoarcs, pseudocircles and Lakes of Wada continuum appear
naturally in a typical map from $S^2$ into $I$.

Our second topic is the ``worst case behavior'' of the level sets
structure of smooth functions. This will involve geometric measure
theory and descriptive set theory. We will prove recent results
which describe how to ``parametrize'' Hausdorff dimension of
analytic sets by smooth functions. As a simple corollary to these
results, it will follow that there is a $C^{\infty}$ function $f:I
\rightarrow I$ such that the set of points where the level sets of
$f$ is uncountable is large in terms of measure and complicated
descriptive set theoretically. This is counter intuitive to the
belief held by some that ``\textit{$C^{\infty}$ functions are more or less
real analytic.}'' We will discuss many interesting and open
question.

The indicatrix of a function is the function
assigning to a point in the range of a function the cardinality of its
level set.
The behaviour of indicatrices of Lebesgue measurable functions
and Borel measurable functions will be described. This seems to be
a topic of still many open and interesting questions.

Here is the plan of the second half of the workshop:

\begin{enumerate}
\item Baire's category theorem and Banach's proof of the existence of a
nowhere differentiable function.

\item Pompeiu derivarives and K\"opcke (i.e. everywhere differentiable and
nowhere monotone) functions.
Weil's proof of the existence of K\"opcke functions using Baire's category
method.

\item Level sets of `typical' continuous functions in one dimension - the
theorem of Bruckner and Garg.

\item Level sets of `typical' $C^n$ functions in one dimension and of typical
continuous functions in two dimensions
- some newer developments.

\item An introduction to Hausdorff measures. The description of collections
of perfect and uncountable level sets for $C^n$ functions.

\item Characterizations of indicatrices of Lebesgue and Borel measurable
functions.
\end{enumerate}
(The results in sections 4,5,6 are due to D'Aniello, Buczolich,
Komisarski, Milewski, Michalewski, Ryll-Nardzewski and to the speakers.)

\subsubsection*{Organizational and additional details}
See
\url{http://www.im.pwr.wroc.pl/\textasciitilde cichon/Bedlewo/}
For additional information contact \texttt{morayne@im.pwr.wroc.pl}

\nby{Jacek Cicho\'n}

\section{Research announcements}

\nAMSPaper{http://www.ams.org/journal-getitem?pii=S0002-9939-05-07861-5}
{Cardinal restrictions on some homogeneous compacta}
{Istvan Juhasz, Peter Nyikos, and Zoltan Szentmiklossy}
{We give restrictions on the cardinality of compact
Hausdorff homogeneous spaces that do not use other cardinal
invariants, but rather covering and separation properties.  In
particular, we show that it is consistent that every hereditarily
normal homogeneous compactum is of cardinality $\mathfrak{c}$.  We introduce
property wD($\kappa $), intermediate between the properties of
being weakly $\kappa $-collectionwise Hausdorff and strongly
$\kappa $-collectionwise Hausdorff, and show that if $X$ is a compact
Hausdorff homogeneous space in which every subspace has property
wD($\aleph _{1}$), then $X$ is countably tight and hence of cardinality
$\le 2^{\mathfrak{c}}$. As a corollary, it is consistent that such a space $X$
is first countable and hence of cardinality $\mathfrak{c}$. A number of
related results are shown and open problems presented.
}

\nArxPaper{math.GN/0504003}
{Filters: Topological congruence relations on groups}
{G\'abor Luk\'acs}
{A filter $\mathcal{F}$ on a group $G$ is a {\em $T$-filter} if there is a
Hausdorff group topology $\tau$ on $G$ such that $\mathcal{F}\stackrel \tau
\longrightarrow 0$. This notion can be specialized for sequences, in which case
we say that $\{a_n\}$ is a {\em $T$-sequence}. In this paper, $T$-filters and
$T$-sequences are studied. We characterize $T$-filters in non-abelian groups,
show that certain filters can be interpreted as topological extensions of the
notion of kernel (i.e., normal subgroup, congruence relation), and provide
several sufficient conditions for a sequence in an abelian group to be a
$T$-sequence. As an application, special sequences in the Pr\"ufer groups
$\mathbb{Z}(p^\infty)$ are investigated. We prove that for $p\neq 2$, there is
a Hausdorff group topology $\tau$ on $\mathbb{Z}(p^\infty)$ that is neither
maximally nor minimally almost periodic--in other words, the von Neumann
radical $\mathbf{n}(\mathbb{Z}(p^\infty),\tau)$ is a non-trivial finite
subgroup. In particular, $\mathbf{n}(\mathbf{n}(\mathbb{Z}(p^\infty),\tau))
\subsetneq \mathbf{n}(\mathbb{Z}(p^\infty),\tau)$.
}

\nArxPaper{math.GN/0504214}
{Inverse Limits and Function Algebras}
{Joan E.\ Hart and Kenneth Kunen}
{Assuming Jensen's principle diamond, there is a compact Hausdorff space $X$
which is hereditarily Lindel\"of, hereditarily separable, and connected, such
that no closed subspace of $X$ is both perfect and totally disconnected. The
Proper Forcing Axiom implies that there is no such space. The diamond example
also fails to satisfy the CSWP (the complex version of the Stone-Weierstrass
Theorem). This space cannot contain the two earlier examples of failure of the
CSWP, which were totally disconnected -- specifically, the Cantor set (W.\ Rudin)
and $\beta\N$ (Hoffman and Singer).}

\nArxPaper{math.LO/0504199}
{Ultrafilters and partial products of infinite cyclic groups}
{Andreas Blass and Saharon Shelah}
{We consider, for infinite cardinals $\kappa$ and $\alpha\le\kappa^+$, the group
$\Pi(\kappa,<\alpha)$ of sequences of integers, of length $\kappa$, with non-zero
entries in fewer than alpha positions. Our main result tells when
$\Pi(\kappa,<\alpha)$ can be embedded in $\Pi(\lambda,<\beta)$. The proof involves some
set-theoretic results, one about families of finite sets and one about families
of ultrafilters.}

\nArxPaper{math.LO/0504200}
{More on regular reduced products}
{Juliette Kennedy and Saharon Shelah}
{The authors show, by means of a finitary version $\square^{fin}_{\lambda,D}$ of
the combinatorial principle $\square^{b^*}_{\lambda}$, the consistency of the
failure, relative to the consistency of supercompact cardinals, of the
following: For all regular filters $D$ on a cardinal $\lambda$, if $M_i$ and $N_i$ are
elementarily equivalent models of a language of size $\le\lambda$, then the second
player has a winning strategy in the Ehrenfeucht-Fraisse game of length
$\lambda^+$ on $\prod_i M_i/D$ and $\prod_i N_i/D$. If in addition $2^{\lambda}= \lambda^+$
and $i<\lambda$ implies $|M_i|+|N_i| \le \lambda^+$, this means that the ultrapowers
are isomorphic.
}

\nArxPaper{math.LO/0504201}
{Consistency of ``the ideal of null restricted to some $A$ is $\kappa$-complete
not $\kappa^+$-complete, $\kappa$ weakly inaccessible and $\cov(\M)=\aleph_1$''}
{Saharon Shelah}
{In this note we answer the following question of Grinblat:
Is it consistent that for some set $A$, $\cov(null\restriction A)=\lambda$ is a
weakly inaccessible cardinal (so $A$ is not null) while $\cov(meager)$ is
small, say it is $\aleph_1$.
}

\nArxPaper{math.GN/0504325}
{On removing one point from a compact space}
{Gady Kozma}
{If $B$ is a compact space and $B\sm\{pt\}$ is Lindel\"of then $B^k\sm\{pt\}$ is star-Linedl\"of
for every cardinality $k$. If $B\sm\{pt\}$ is compact then $B^k\sm\{pt\}$ is discretely
star-Lindel\"of. In particular, this gives new examples of Tychonoff discretely
star-Lindel\"of spaces with unlimited extent.}

\section{Problem of the Issue}

Since this issue is sent out earlier than planned, it will contain no
\emph{Problem of the Issue} this time.
The next section contains many problems which are still open,
so that the reader may wish to solve any of them instead.

\section{Problems from earlier issues}
In this section we list the still open problems among
the past problems posed in the \SPMBul{}
(in the section \emph{Problem of the month/issue}).
For definitions, motivation and related results, consult the
corresponding issue.

For conciseness, we make the convention that
all spaces in question are
zero-dimentional, separable metrizble spaces.

\begin{issue}
Is $\binom{\Omega}{\Gamma}=\binom{\Omega}{\Tau}$?
\end{issue}

\begin{issue}
Is $\ufin(\Gamma,\Omega)=\sfin(\Gamma,\Omega)$?
And if not, does $\ufin(\Gamma,\Gamma)$ imply
$\sfin(\Gamma,\Omega)$?
\end{issue}

\stepcounter{issue}

\begin{issue}
Does $\sone(\Omega,\Tau)$ imply $\ufin(\Gamma,\Gamma)$?
\end{issue}

\begin{issue}
Is $\p=\p^*$? (See the definition of $\p^*$ in that issue.)
\end{issue}

\begin{issue}
Does there exist (in ZFC) an uncountable set satisfying $\sone(\BG,\B)$?
\end{issue}

\stepcounter{issue}

\begin{issue}
Does $X \nin \NON(\M)$ and $Y\nin\mathsf{D}$ imply that
$X\cup Y\nin \COF(\M)$?
\end{issue}

\begin{issue}
Is $\split(\Lambda,\Lambda)$ preserved under taking finite unions?
\end{issue}
\begin{proof}[Partial solution]
Consistently yes (Zdomsky). Is it ``No'' under CH?
\end{proof}

\begin{issue}
Is $\cov(\M)=\fo$? (See the definition of $\fo$ in that issue.)
\end{issue}

\begin{issue}
Does $\sone(\Gamma,\Gamma)$ always contain an element of cardinality $\b$?
\end{issue}

\begin{issue}
Could there be a Baire metric space $M$ of weight $\aleph_1$ and a partition
$\mathcal{U}$ of $M$ into $\aleph_1$ meager sets where for each ${\mathcal U}'\subset\mathcal U$,
$\bigcup {\mathcal U}'$ has the Baire property in $M$?
\end{issue}

\general\end{document}